\documentclass[journal,twoside,web]{IEEEtran}
\usepackage{generic}
\usepackage{cite}
\usepackage{amsmath,amssymb,amsfonts}
\usepackage{algorithmic}
\usepackage{graphicx}
\usepackage{algorithm,algorithmic}
\usepackage{hyperref}
\hypersetup{hidelinks=true}
\usepackage{textcomp}

\newtheorem{lemma}{Lemma}
\newtheorem{assumption}{Assumption}
\newtheorem{remark}{Remark}
\newtheorem{theorem}{Theorem}

\def\BibTeX{{\rm B\kern-.05em{\sc i\kern-.025em b}\kern-.08em
    T\kern-.1667em\lower.7ex\hbox{E}\kern-.125emX}}
\begin{document}
	
\title{Actor–Critic Learning for Risk-Constrained Linear Quadratic Regulation}

\author{Weijian Li and Andreas A. Malikopoulos 
\thanks{This research was supported in part by NSF under Grants CNS-2401007, CMMI-2219761, IIS-2415478, and in part by MathWorks.}
\thanks{Weijian Li and Andreas A. Malikopoulos are with the Department of Civil and Environmental Engineering, Cornell University, Ithaca, NY, USA (e-mails: wl779@cornell.edu, amaliko@cornell.edu).}
}

\maketitle

\begin{abstract}
In this paper, we investigate the infinite-horizon risk-constrained linear quadratic regulator problem (RC-LQR), which augments the classical LQR formulation with a statistical constraint on the variability of the system state to incorporate risk awareness, a key requirement in safety-critical control applications.
We propose an actor–critic learning algorithm that jointly performs policy evaluation and policy improvement in a model-free and online manner. The RC-LQR problem is first reformulated as a max–min optimization problem, from which we develop a multi–time–scale stochastic approximation scheme. The critic employs temporal-difference learning to estimate the action–value function, the actor updates the policy parameters via a policy gradient step, and the dual variable is adapted through gradient ascent to enforce the risk constraint.
\end{abstract}

\begin{IEEEkeywords}
Risk-constrained LQR,
actor-critic algorithm,
stochastic approximation,
reinforcement learning.
\end{IEEEkeywords}

\section{Introduction}
\label{sec:introduction}

Reinforcement learning (RL), aiming to learn optimal policies directly from data in the absence of explicit model knowledge, has gained significant attention in recent years. 
Its versatility has enabled progress across a wide range of sequential decision-making problems, including poker games~\cite{brown2019superhuman}, robotic control~\cite{kober2013reinforcement}, and cyber–physical systems~\cite{leong2020deep}. 
Recent studies have investigated the connection between learning and control by introducing information–state representations that decouple estimation and control~\cite{Malikopoulos2022a, Malikopoulos2024}, and by developing integrated learning–and–control frameworks for systems with unknown dynamics~\cite{kounatidis2025combined}. 
Within the control community, there has been an increased interest in establishing a rigorous theoretical foundation for RL in linear systems—most notably through the \emph{linear quadratic regulator} (LQR) problem, which seeks to regulate an unknown linear system to minimize a quadratic cost function~\cite{recht2019tour, jiang2017robust}.

The LQR problem provides a mathematically tractable framework for studying the interplay between learning and control and serves as a canonical benchmark for analyzing the stability, optimality, and convergence of data-driven control algorithms. 
However, the classical formulation is \emph{risk-neutral}, as performance is evaluated solely through the expected cost, overlooking low-probability but high-impact disturbances that can cause unsafe behaviors in safety-critical systems~\cite{whittle1981risk}. 
This shortcoming has motivated the development of \emph{risk-aware} formulations that explicitly regulate the variability of the closed-loop behavior. 
Early work introduced an exponential cost transformation to capture risk sensitivity~\cite{whittle1981risk, pan1996model, cui2024robust}, though such problems can become ill-posed under general noise distributions. 
Conditional value-at-risk formulations~\cite{chapman2019risk, chapman2021risk} offer a more coherent treatment of tail events but often rely on tractable approximations. 
The \emph{risk-constrained LQR (RC-LQR)} framework~\cite{tsiamis2020risk} augments the standard LQR objective with a constraint on the predictive variance of the state penalty, yielding a closed-form optimal policy that accommodates arbitrary noise models. 
Recent extensions established global convergence guarantees for policy optimization methods~\cite{zhao2021infinite, zhao2023global}.

Advances in reinforcement learning (RL) have been largely driven by \emph{policy optimization} methods, such as policy gradient and actor–critic algorithms~\cite{sutton1999policy, konda1999actor, bhatnagar2009natural, schulman2015trust}, which directly parameterize and optimize control policies. 
These approaches efficiently handle high-dimensional, continuous spaces and are simple to implement. 
Theoretical analyses have also been developed for linear control systems~\cite{fazel2018global, hu2023toward, fatkhullin2021optimizing, duan2023optimization}, while primal–dual policy gradient algorithms have been proposed for the RC-LQR problem~\cite{zhao2023global}. 
Among these methods, the \emph{actor–critic} framework, which performs policy evaluation and improvement simultaneously, has received extensive attention~\cite{bhatnagar2009natural, zhang2018fully, zeng2024two, yang2019provably}. 
Compared with sample-based policy gradient schemes, actor–critic algorithms typically converge faster due to variance reduction in gradient estimation. 
However, their extension to \emph{constrained} RL problems remains relatively unexplored.

In this paper, we develop an actor–critic algorithm for the RC-LQR problem and characterize its theoretical properties. 
The main contributions of this work are summarized as follows:
\textbf{(1)} We study a \emph{model-free} LQR problem subject to a risk constraint. 
Unlike the classical LQR formulation~\cite{jiang2017robust, fazel2018global, yang2019provably}, the proposed setting minimizes a quadratic cost and also regulates the statistical variability of the system state. 
While similar formulations appear in~\cite{tsiamis2020risk, zhao2021infinite, patel2024risk}, our focus is on the model-free regime.
\textbf{(2)} We propose a \emph{multi–time–scale actor–critic algorithm} under a \emph{primal–dual} framework. 
The critic employs temporal-difference learning, the actor performs policy-gradient updates, and the dual variable evolves via gradient ascent. 
The proposed scheme generalizes two–time–scale actor–critic algorithms~\cite{bhatnagar2009natural, zhang2018fully, zeng2024two} to constrained RL and, unlike the sample-based policy gradient method of~\cite{zhao2023global}, integrates policy evaluation and improvement in an online and model-free manner.
\textbf{(3)} Using stochastic approximation theory, we establish \emph{almost sure asymptotic convergence} of all iterates. 
This result extends classical actor–critic convergence guarantees~\cite{bhatnagar2009natural, zhang2018fully, zeng2024two} by incorporating a dual update in the learning loop and strengthens prior results such as~\cite{zhao2023global}, which established convergence only with high probability.

The remainder of the paper is structured as follows. 
In Section~\uppercase\expandafter{\romannumeral2}, we review the necessary preliminaries that form the basis of our analysis. 
In Section~\uppercase\expandafter{\romannumeral3}, we present the formulation of the RC-LQR problem. 
In Section~\uppercase\expandafter{\romannumeral4}, we develop the proposed actor–critic algorithm under a primal–dual framework, and in Section~\uppercase\expandafter{\romannumeral5}, we establish its asymptotic convergence using stochastic approximation theory. 
In Section~\uppercase\expandafter{\romannumeral6}, we provide simulation results that illustrate the performance of the algorithm, and in Section~\uppercase\expandafter{\romannumeral7}, we draw concluding remarks.

\textbf{Notation:}
Let $\mathbb{N} := \{1,2,\dots\}$ denote the set of positive integers, and fix dimensions $n,m \in \mathbb{N}$. 
Let $\mathbb{R}_{+} := [0,\infty)$, $\mathbb{R}^{n}$ be the set of $n$-dimensional real column vectors, $\mathbb{R}^{m \times n}$ be the set of $m \times n$ real matrices, and $I_{n}$ the $n \times n$ identity. 
Let $\mathbb{S}^{n} \subset \mathbb{R}^{n \times n}$ denote the set of real $n \times n$ symmetric matrices, and write $\mathbb{S}_{++}^{n}$ ($\mathbb{S}_{+}^{n}$) for the positive (semi-)definite cones in $\mathbb{S}^{n}$. 
We use $X \succ (\succeq)\,0$ to indicate $X \in \mathbb{S}_{++}^{n}$ ($\mathbb{S}_{+}^{n}$). The transpose and trace operators are $(\cdot)^{\top}$ and $\mathrm{tr}(\cdot)$, respectively. 
For $A,B \in \mathbb{R}^{m \times n}$, the Frobenius inner product is $\langle A,B\rangle_{F} := \mathrm{tr}(A^{\top}B)$, with Frobenius norm $\|A\|_{F} := \sqrt{\langle A,A\rangle_{F}}$. 
For $x \in \mathbb{R}^{n}$, $\|x\|$ denotes the Euclidean ($\ell_{2}$) norm. 
For $X \in \mathbb{S}^{n}$, $\mathrm{svec}(X)$ denotes the vectorization of the upper triangular part of $X$, with off-diagonal entries weighted by $\sqrt{2}$. 
Given a nonempty closed convex set $\mathcal{X} \subset \mathbb{R}^{n}$, the Euclidean projection is $
\Gamma_{\mathcal{X}}(y) := \arg\min_{x \in \mathcal{X}} \|y - x\|,$ for all $y \in \mathbb{R}^{n}.$

\section{Preliminaries}
\label{sec:preliminary}

In this section, we review fundamental results from stochastic approximation theory that will be used in the convergence analysis. 

Consider a sequence $\{x_t\} \subset \mathbb{R}^n$ generated by
\begin{equation}
\label{SA:alg}
x_{t+1} = x_t + \gamma_t [h(x_t) + M_{t+1}], \qquad t \ge 0,~ x_0 \in \mathbb{R}^n,
\end{equation}
where $h(\cdot)$ is a deterministic map, $\{M_t\}$ is a noise sequence, and $\{\gamma_t\}$ denotes the step size.

\begin{assumption}
\label{ass:SA}
The recursion~\eqref{SA:alg} satisfies:
\begin{enumerate}
\item $h:\mathbb{R}^n \!\to\! \mathbb{R}^n$ is Lipschitz continuous.
\item The step sizes obey $\sum_t \gamma_t = \infty$ and $\sum_t \gamma_t^2 < \infty$.
\item $\{M_t\}$ is a martingale difference sequence with $M_0=0$ and $\mathbb{E}[M_{t+1}\mid x_s,M_s,s\!\le\! t]=0$, and there exists $K>0$ such that
\[
\mathbb{E}[\|M_{t+1}\|^2 \mid x_s,M_s,s\!\le\! t] \le K(1+\|x_t\|^2).
\]
\end{enumerate}
\end{assumption}

The asymptotic behavior of~\eqref{SA:alg} is characterized by the limiting ordinary differential equation (ODE)
\begin{equation}
\label{SA:ODE}
\dot{x} = h(x).
\end{equation}
If~\eqref{SA:ODE} admits a unique globally asymptotically stable equilibrium $x^*$, the following results hold~\cite{borkar2008stochastic}.

\begin{lemma}
\label{lem:SA:con}
Under Assumption~\ref{ass:SA}, if $\sup_t \|x_t\| < \infty$ a.s., then $x_t \to x^*$ a.s.
\end{lemma}

\begin{lemma}
\label{lem:SA:bound}
Let $h_{\infty}(x) := \lim_{r\to \infty} {h(rx)}/{r}$. Under Assumption \ref{ass:SA}, suppose $h_{\infty}$
exists uniformly on compact sets for some $h_{\infty} \in \mathcal C(\mathbb R^n)$. If the ODE $\dot y = h_{\infty}(y)$ has origin as its unique globally asymptotically stable equilibrium, then 
$\sup_t \Vert x_t\Vert < \infty$ a.s.
\end{lemma}

Next, we recall the Kushner–Clark lemma~\cite{kushner2003stochastic}, which addresses projected stochastic approximation.

Let $\mathcal{X} \subset \mathbb{R}^n$ be a compact convex set and define
\[
\widehat{\Gamma}_{\mathcal{X}}[h(x)] := \lim_{0 < \eta \to 0} 
\frac{\Gamma_{\mathcal{X}}[x + \eta h(x)] - x}{\eta}, \qquad x \in \mathcal{X}.
\]
Consider the recursion
\begin{equation}
\label{SA:alg2}
x_{t+1} = \Gamma_{\mathcal{X}}\!\big[x_t + \gamma_t (h(x_t) + \xi_t + \beta_t)\big],
\end{equation}
whose limiting ODE is
\begin{equation}
\label{SA:ODE2}
\dot{x} = \widehat{\Gamma}_{\mathcal{X}}[h(x)].
\end{equation}

\begin{assumption}
\label{ass:SA2}
The recursion~\eqref{SA:alg2} satisfies:
\begin{enumerate}
\item $h:\mathbb{R}^n \!\to\! \mathbb{R}^n$ is continuous.
\item $\{\beta_t\}$ is bounded and $\beta_t \to 0$ a.s.
\item $\gamma_t \to 0$, $\sum_t \gamma_t = \infty$.
\item For any $\varepsilon > 0$,
\[
\lim_t \mathbb{P}\!\left(\sum_{n \ge t} 
\Big\|\!\sum_{s=t}^n \gamma_s \xi_s\!\Big\|^2 \ge \varepsilon \right) = 0.
\]
\end{enumerate}
\end{assumption}

If the ODE~\eqref{SA:ODE2} admits a compact asymptotically stable equilibrium set $\mathcal{K}^*$, the following result holds.

\begin{lemma}
\label{lem:SA:proj}
Under Assumption~\ref{ass:SA2}, the iterates $x_t$ in~\eqref{SA:alg2} satisfy $x_t \to \mathcal{K}^*$ a.s. as $t \to \infty$.
\end{lemma}

\section{Problem Formulation}

Consider the discrete-time linear time-invariant (LTI) stochastic system
\begin{equation}
\label{LTI:sys}
x_{t+1} = A x_t + B u_t + w_t,
\end{equation}
where $x_t \in \mathbb{R}^n$ denotes the state, $u_t \in \mathbb{R}^m$ the control input, and $\{w_t\}$ an independent and identically distributed (i.i.d.) noise sequence. 
The system matrices $A \in \mathbb{R}^{n \times n}$ and $B \in \mathbb{R}^{n \times m}$ are unknown but constant. 
The infinite-horizon linear quadratic regulator (LQR) problem seeks a control policy $u := \{u_t\}_{t \ge 0}$ that minimizes the long-run average cost
\begin{equation}
\label{cost}
\mathcal{J}(u) := \lim_{T \to \infty} \frac{1}{T}\,
\mathbb{E} \!\left[\sum_{t=0}^{T-1} \!\left(x_t^{\top} Q x_t + u_t^{\top} R u_t \right)\!\right],
\end{equation}
where $Q \in \mathbb{S}_{+}^{n}$ and $R \in \mathbb{S}_{++}^{m}$ are the state and input weighting matrices, respectively.

The LQR problem is \emph{risk-neutral}, as it optimizes only the expected performance. 
Consequently, the state may grow arbitrarily large by extreme noise realizations, particularly when $\{w_t\}$ follows a heavy-tailed distribution. 
To incorporate risk sensitivity, following~\cite{tsiamis2020risk, zhao2023global, patel2024risk}, we consider the \emph{risk-constrained LQR (RC-LQR)} problem formulated as
\begin{equation}
\begin{aligned}
\label{formulation}
\min_{u}~&\lim_{T \to \infty}\frac{1}{T}\,
\mathbb{E}\!\left[\sum_{t=0}^{T-1}\!\left(x_t^{\top} Q x_t + u_t^{\top} R u_t\right)\!\right] \\
{\rm s.t.}~&
x_{t+1} = A x_t + B u_t + w_t, \\
&\lim_{T \to \infty}\frac{1}{T}\,
\mathbb{E}\!\left[\sum_{t=0}^{T-1}\!
\Big(x_t^{\top} Q x_t - 
\mathbb{E}[x_t^{\top} Q x_t \mid \mathcal{F}_{t-1}]\Big)^{\!2}\right] 
\le \iota,
\end{aligned}
\end{equation}
where $\iota > 0$ is a user-defined risk-tolerance parameter, and 
$\mathcal{F}_t := \sigma\{x_s, u_s : s \le t\}$ denotes the $\sigma$-algebra generated by all the observations up to time $t$.

To ensure the well-posedness of~\eqref{formulation}, we impose the following standing assumptions.

\begin{assumption}
\label{ass:mat}
The pair $(A,B)$ is stabilizable, $Q \in \mathbb{S}_{+}^{n}$, $R \in \mathbb{S}_{++}^{m}$, and the pair $(A, Q^{1/2})$ is detectable.
\end{assumption}

\begin{assumption}
\label{ass:noise}
The noise sequence $\{w_t\}$ satisfies $\mathbb{E}\|w_t\|^4 < \infty$.
\end{assumption}

\begin{assumption}
\label{ass:slater}
The Slater condition holds; namely, there exists a feasible control policy that strictly satisfies the risk constraint in~\eqref{formulation}.
\end{assumption}

\begin{remark}
Unlike the standard LQR problem, \eqref{formulation} incorporates the predictive variance of the state penalty, thereby enhancing risk awareness. 
Decreasing the risk-tolerance parameter~$\iota$ tightens the constraint, providing a systematic trade-off between average performance and risk sensitivity.
\end{remark}

Let $\bar{w} := \mathbb{E}[w_t]$ and 
$W := \mathbb{E}[(w_t - \bar{w})(w_t - \bar{w})^{\top}]$ 
denote the mean and covariance of the noise, respectively, and assume $W \succ 0$. 
Define the higher-order moments
\begin{equation*}
\begin{aligned}
M_3 &:= \mathbb{E}\!\left[(w_t - \bar{w})(w_t - \bar{w})^{\top} 
Q (w_t - \bar{w})\right], \\
m_4 &:= \mathbb{E}\!\left[\big((w_t - \bar{w})^{\top} Q (w_t - \bar{w}) 
- \mathrm{tr}(WQ)\big)^{2}\right].
\end{aligned}
\end{equation*}

By~\cite[Theorem~1]{zhao2021infinite}, the optimal policy of~\eqref{formulation} admits a time-invariant affine structure, i.e.,
$u_t = -K^{*} x_t + b^{*},$
where $K^{*} \in \mathbb{R}^{m \times n}$ and $b^{*} \in \mathbb{R}^{m}$. 
Hence, it is without loss of optimality to restrict attention to affine policies of the form $u_t = -Kx_t + b$. 
Let $X := [K,\, b]$ denote the decision variable of~\eqref{formulation}, and $X^{*} := [K^{*},\, b^{*}]$ its optimal value. 
Following the argument in~\cite[Proposition~1]{tsiamis2020risk}, problem~\eqref{formulation} can be equivalently reformulated as
\begin{equation}
\begin{aligned}
\label{reform}
\min_{X}~& \mathcal{J}(X) \\
{\rm s.t.}~~& x_{t+1} = A x_t + B u_t + w_t, \\
& u_t = -Kx_t + b, \\
& \mathcal{J}_{c}(X) \le \bar{\iota},
\end{aligned}
\end{equation}
where $\bar{\iota} := \iota - m_4 + 4\,\mathrm{tr}[(WQ)^2]$, and
\begin{equation*}
\mathcal{J}(X)
:= \lim_{T \to \infty} \frac{1}{T}\,
\mathbb{E}\!\left[\sum_{t=0}^{T-1}\!\left(x_t^{\top} Q x_t + u_t^{\top} R u_t\right)\!\right],
\end{equation*}
\begin{equation*}
\mathcal{J}_{c}(X)
:= \lim_{T \to \infty} \frac{1}{T}\,
\mathbb{E}\!\left[\sum_{t=0}^{T-1}\!
\left(4\,x_t^{\top} QWQ x_t + 4\,x_t^{\top} Q M_3\right)\!\right].
\end{equation*}

\begin{remark}
The RC-LQR problem was originally introduced in~\cite{tsiamis2020risk, zhao2021infinite} and later extended to mean-field linear–quadratic systems in~\cite{patel2024risk}. 
Subsequent work~\cite{zhao2023global} developed primal–dual policy gradient algorithms for this formulation. 
In contrast, the present paper proposes an \emph{actor–critic algorithm} that solves the RC-LQR problem in an online, model-free setting. 
The method jointly performs policy evaluation and improvement within a variance-reduction framework, achieving faster convergence than standard policy gradient approaches.
\end{remark}

\section{Actor-Critic Algorithm for RC-LQR}

In this section, we analyze the properties of the state- and action-value functions and establish a policy gradient representation that links policy gradients to the action-value function. Building on these results, we develop a multi–time–scale actor–critic algorithm under a primal–dual framework.

\subsection{Supporting Results}

Let the set of stabilizing policies for~\eqref{LTI:sys} be
\[
\mathcal{S} := \{[K,b] \mid \rho(A - BK) < 1,\; K \in \mathbb{R}^{m \times n},\; b \in \mathbb{R}^{m}\},
\]
where $\rho(\cdot)$ denotes the spectral radius. 

The Lagrangian function associated with~\eqref{reform} is
\begin{equation}
\begin{aligned}
\label{lagrangian}
\mathcal{L}(X,\mu)
&= \mathcal{J}(X) + \mu\big(\mathcal{J}_{c}(X) - \bar{\iota}\big) \\
&= \lim_{T \to \infty} \frac{1}{T}\,
\mathbb{E}\!\left[\sum_{t=0}^{T-1} c_{\mu}(x_t,u_t)\right],
\end{aligned}
\end{equation}
where $\mu \in \mathbb{R}_{+}$ is the dual variable,
\[
Q_{\mu} := Q + 4\mu QWQ,\qquad S := 2\mu Q M_{3},
\]
and
\[
c_{\mu}(x_t,u_t)
:= x_t^{\top} Q_{\mu} x_t + 2 x_t^{\top} S + u_t^{\top} R u_t - \mu \bar{\iota}.
\]

By~\cite[Theorem~2]{zhao2023global}, there is no duality gap between~\eqref{reform} and its dual problem. Hence, the optimal solution can be obtained by solving
\begin{equation}
\label{duality}
\max_{\mu \ge 0}\; \min_{X \in \mathcal{S}}\, \mathcal{L}(X,\mu).
\end{equation}

Throughout the sequel, for any $\mu \ge 0$, we consider the family of linear–Gaussian policies
\begin{equation}
\label{policy}
\Pi^{\mu}
:= \{\pi_{X}^{\mu}(\cdot \mid x)
= \mathcal{N}(-Kx + b,\, \sigma^{2} I_{m}) :
K \in \mathbb{R}^{m \times n},\, b \in \mathbb{R}^{m}\},
\end{equation}
where, for each $t \ge 0$, the control input is $
u_t = -Kx_t + b + \sigma \eta_t,~ \eta_t \sim \mathcal{N}(0, I_m).
$
The optimal deterministic policy $u_t = -K^{*}x_t + b^{*}$ is recovered by taking $\sigma = 0$. 
The injected Gaussian noise promotes exploration during learning while preserving convergence to the optimal affine policy.

For any policy $\pi^\mu_X \in \Pi^\mu$, it follows from (\ref{LTI:sys}) that
\begin{equation}
\label{LTI:reform}
x_{t+1} = (A - BK)x_t + Bb + \zeta_t,
\end{equation}
where $\zeta_t = w_t + \sigma B \eta_t$.
Let $\bar \zeta$ and $\Psi_\zeta$ be the mean and covariance of $\zeta_t$. Since $w_t$ and $\eta_t$ are independent, we have $\bar \zeta = \bar w$, and $\Psi_\zeta = W + \sigma^2 BB^\top$.
It is known that, for any $X \in \mathcal S$, the Markov chain (\ref{LTI:reform}) admits a stationary distribution, denoted by $\nu_X$, whose mean $\bar x_X$ and covariance $\Sigma_{K}$ satisfy
\begin{equation}
\begin{aligned}
\bar x_{X} &= (A - BK) \bar x_X + Bb + \bar w, \\
\Sigma_K &= \Psi_\zeta + (A - BK) \Sigma_K (A - BK)^\top.
\end{aligned}
\end{equation}

We define the state- and action-value functions associated with the Lagrangian (\ref{lagrangian}) as 
\begin{equation}
\begin{aligned}
V^{\mu}_{X}(x) &= \mathbb E \left[\sum_{t = 0}^\infty \big(c_\mu(x_t, u_t) - \mathcal L(X, \mu)\big) \mid x_0 = x\right],\\
Q^{\mu}_X(x, u) &= \mathbb E \left[\sum_{t = 0}^\infty \big(c_\mu(x_t, u_t) \!-\! \mathcal L(X, \mu)\big) \!\mid \! x_0 \!=\! x, u_0 \!=\! u\right], 
\end{aligned}
\end{equation}
where $\mathbb E [\cdot]$ takes the expectation under a fixed policy $X \in \mathcal S$ and a fixed multiplier $\mu \in \mathbb R_+$.

The following lemma indicates that both the state- and action-value functions are quadratic.
\begin{lemma}
\label{lem:value:function}
For any $X \in \mathcal S$, and $\mu \in \mathbb R_+$, the following statements hold.
\begin{enumerate}
\item
\makebox[\linewidth][c]{$
V^{\mu}_X(x) = x^\top P_{K,\mu} x + g^\top_{X,\mu} x + z^1_{X,\mu},$}
where $P_{K,\mu}$ is the solution to a Lyapunov function as
\begin{equation*}
\begin{aligned}
P_{K, \mu} &= Q_\mu + K^\top R K + (A - BK)^\top P_{K,\mu} (A - BK), \\
g_{X,\mu} &= 2[I - (A - BK)]^{-\top}\big(S - K^\top R b \\
&+ (A - BK)^\top P_{K,\mu} B b + (A - BK)^\top P_{K,\mu} \bar w\big),
\end{aligned}
\end{equation*}
and $z^{1}_{X,\mu}$ is a constant independent of $x$.

\item 
\begin{equation}
\label{Q:form}
\begin{aligned}
Q^{\mu}_{X}&(x, u) \!=\! 
\begin{bmatrix}
x \\
u
\end{bmatrix}^\top \!\!\Upsilon_{K,\mu} \!\!
\begin{bmatrix}
x \\
u
\end{bmatrix}
\!+\! 2\! \begin{bmatrix}
p_{X,\mu} \\
q_{X,\mu}
\end{bmatrix}^\top \!\!
\begin{bmatrix}
x \\
u
\end{bmatrix} \!+\! z^2_{X,\mu},
\end{aligned}
\end{equation}
where
\begin{equation*}
\begin{aligned}
\Upsilon_{K,\mu} \!\!=\!\! 
\begin{bmatrix}
\Upsilon_{K,\mu}^{11}~\Upsilon_{K,\mu}^{12} \\
\Upsilon_{K,\mu}^{21}~\Upsilon_{K,\mu}^{22}
\end{bmatrix}
\!\!=\!\! \begin{bmatrix}
Q_\mu \!+\! A^\top P_{K,\mu} A  
\!\!\!\!\!\!&A^\top P_{K,\mu} B \\
B^\top P_{K,\mu} A		
\!\!\!\!\!\!&R \!+\! B^T P_{K,\mu} B
\end{bmatrix}\!\!, 
\end{aligned}
\end{equation*}
$p_{X,\mu} =  S + A^\top P_{K,\mu} \bar w + \frac 12 A^\top g_{X,\mu}$, 
$q_{X,\mu} = B^\top P_{K,\mu} \bar w + \frac 12 B^\top g_{X,\mu}$,
and $z^2_{X, \mu}$ is a constant independent of $(x, u)$.
\end{enumerate}
\end{lemma}

\emph{Proof:}
Substituting $u_t = -Kx_t + b + \sigma \eta_t$ into $c_{\mu}(x_t, u_t)$ yields
\begin{equation*}
\begin{aligned}
\mathbb E & \left[c_{\mu}(x_t,  u_t) \mid x_t \right] 
= x_t^\top (Q_\mu + K^\top R K)x_t \\
&+ 2x_t^\top(S - K^\top R b) + \sigma^2 {\rm tr}(R) + b^T R b - \mu \bar \iota.
\end{aligned}
\end{equation*}
It follows that
\begin{equation*}
\label{L:form}
\begin{aligned}
\mathcal L(&X,\mu) =~ \lim_{T \to \infty} \mathbb E \left[\frac 1 T \sum_{t = 0}^{T-1} \mathbb E \big[c_{\mu}(x_t, u_t) \mid x_t\big]\right] \\
=&~ {\rm tr}\big[(Q_\mu + K^\top R K) \Sigma_K\big]
+ \bar x_X^\top (Q_\mu + K^\top R K)\bar x_X \\
&+ 2 \bar x_X^\top(S - K^\top R b)  + \sigma^2 {\rm tr}(R) + b^T R b - \mu \bar \iota.
\end{aligned}
\end{equation*}

Following the standard LQR analysis in~\cite[Ch.~3]{bertsekas2012dynamic}, 
the value function $V^{\mu}_{X}(x)$ admits a quadratic representation,
i.e., $V^{\mu}_{X}(x) = x^{\top} P_{K,\mu} x + g_{X,\mu}^{\top} x + z^{1}_{X,\mu}$. 
Moreover, $V^{\mu}_{X}$ satisfies the Bellman equation
\begin{equation*}
\begin{aligned}
V^{\mu}_X(x_t) = \mathbb E \big[&c_\mu(x_t, u_t) - \mathcal L(X, \mu) + V^{\mu}_X(x_{t+1}) \mid x_t, \\
&u_t = -Kx_t + b + \sigma \eta_t \big].
\end{aligned}
\end{equation*}
More precisely, it follows that
\begin{equation*}
\begin{aligned}
&x_t^\top P_{K,\mu} x_t + g^{\top}_{X,\mu} x_t + z^1_{X,\mu} = x_t^\top (Q_\mu + K^\top R K)x_t \\
&+ 2x_t^\top(S - K^\top R b) + \sigma^2 {\rm tr}(R) + b^T R b - \mu \bar \iota - \mathcal L(X, \mu) \\
&+ \mathbb E \left\{[(A \!-\! BK) x_t \!+\! Bb + \zeta_t]^\top P_{K,\mu} [(A \!-\! BK) x_t \!+\! Bb \!+\! \zeta_t] \right\} \\
& + \mathbb E \left\{g^\top_{X, \mu} [(A - BK) x_t + Bb + \zeta_t] \right\} + z^1_{X,\mu}.
\end{aligned}
\end{equation*}
As a result, we obtain
\begin{equation*}
\begin{aligned}
P_{K,\mu} &= Q_\mu + K^\top RK + (A-BK)^\top P_{K,\mu} (A-BK), \\
g_{X,\mu} &=2(S - K^\top R b) + 2(A - BK)^\top P_{K,\mu} B b \\
&~~~~+ 2(A - BK)^\top P_{K,\mu} \bar w + (A - BK)^\top g_{X,\mu}.
\end{aligned}
\end{equation*}

Since $\mathbb E_{x \sim \nu_X}\big[V^{\mu}_X(x)\big] = 0$, we have
\begin{equation*}
\begin{aligned}
z^1_{X,\mu} = - {\rm tr}(P_{K,\mu} \Sigma_K) - \bar x_X^\top P_{K,\mu} \bar x_X - g^\top_{X,\mu} \bar x_X.
\end{aligned}
\end{equation*}
Therefore, the first claim follows.

For any $(x_t, u_t)$, we have
\begin{equation}
\begin{aligned}
\label{pf:QV}
Q^{\mu}_{X}(x_t, u_t)
&= c_{\mu}(x_t, u_t) - \mathcal{L}(X, \mu) \\
&\quad + \mathbb{E}\!\left[V^{\mu}_{X}(x_{t+1}) \mid x_t, u_t\right].
\end{aligned}
\end{equation}

Substituting $c_\mu(x_t, u_t)$ and $V^{\mu}_X(x_{t+1})$ into (\ref{pf:QV}), we obtain
\begin{equation*}
\begin{aligned}
Q^{\mu}_X(&x_t, u_t) =  x_t^\top Q_\mu x_t + 2 x_t^\top S + u_t^T R u_t - \mu \bar \iota \\
&+ \mathbb E \big[(Ax_t + Bu_t + w_t)^\top P_{K,\mu}(Ax_t + Bu_t + w_t) \big] \\
&+ \mathbb E \big[g^{\top}_{X,\mu} (Ax_t + Bu_t + w_t)\big] - \mathcal L(X, \mu) \\
=&\begin{bmatrix}
x_t \\
u_t
\end{bmatrix}^\top \!\Upsilon_{K,\mu}\!
\begin{bmatrix}
x_t \\
u_t
\end{bmatrix}
+ 2 \begin{bmatrix}
p_{X,\mu} \\
q_{X,\mu}
\end{bmatrix}^\top \!
\begin{bmatrix}
x_t \\
u_t
\end{bmatrix} 
- \mathcal L(X, \mu) \\
&-\mu \bar \iota + z^1_{X,\mu} + g_{X,\mu}^\top \bar w + {\rm tr}(P_{K,\mu} \bar w \bar w^\top) + {\rm tr}(P_{K,\mu} W),
\end{aligned}
\end{equation*}
where $\Upsilon_{K,\mu}$, $p_{X,\mu}$ and $q_{X,\mu}$ are defined in (\ref{Q:form}).
Hence, the second claim follows, which completes the proof.
$\hfill\square$

The following lemma characterizes the policy gradient of $\mathcal{L}(X,\mu)$ with respect to $X$.

\begin{lemma}
\label{lem:policy_gradient}
For any $X \in \mathcal S$ and $\mu \in \mathbb R_+$, it holds that
\begin{equation}
\label{policy:gradient}
\nabla_X \mathcal L(X, \mu) = 2 H_{X, \mu} \Phi_X, 
\end{equation}
where $H_{X, \mu} \!=\! \big[E_{K,\mu}, \! G_{X,\mu} \big]$, $\Upsilon_{K, \mu}$ and $q_{X,\mu}$ are given in (\ref{Q:form}),
\begin{equation}
\begin{aligned}
\label{grad:Q}
E_{K, \mu} = \Upsilon_{K,\mu}^{22} K - \Upsilon_{K,\mu}^{21},
~G_{X, \mu} = \Upsilon_{K,\mu}^{22} b - q_{X,\mu},
\end{aligned}
\end{equation}
and moreover,
\begin{equation}
\begin{aligned}
\label{Phi:eq}
\Phi_X &= \lim_{T \to \infty} \frac 1T \mathbb E \left\{\sum_{t = 0}^{T-1} \begin{bmatrix}
x_t \\
-1
\end{bmatrix}
\begin{bmatrix}
x_t \\
-1
\end{bmatrix}^\top \right\} \\
&= \begin{bmatrix}
\Sigma_K + \bar x_X \bar x_X^\top  &-\bar x_X \\
-\bar x_X &1
\end{bmatrix} \succ 0. 
\end{aligned}
\end{equation}
\end{lemma}

\emph{Proof:}
By arguments analogous to those in the proof of Lemma~2 in~\cite{zhao2023global}, we obtain
\begin{equation*}
\nabla_X \mathcal L(X, \mu) = 2 \big[E_{K, \mu}, G_{X, \mu} \big] \Phi_X, 
\end{equation*}
where
\begin{equation*}
\begin{aligned}
E_{K, \mu} &= (R + B^\top P_{K,\mu} B)K - B^\top P_{K,\mu} A, \\
G_{X, \mu} &= (R + B^\top P_{K,\mu} B)b + B^\top P_{K,\mu} \bar w + \frac 12 B^\top g_{X,\mu}.
\end{aligned}
\end{equation*}
From~\eqref{Q:form}, we obtain~\eqref{grad:Q}, completing the proof.
$\hfill\square$

\begin{remark}
The stochastic policy class~\eqref{policy} differs from the deterministic counterparts considered in~\cite{zhao2021infinite, zhao2023global} by incorporating Gaussian exploration noise. 
Nevertheless, the state-value function retains its quadratic structure, and the resulting policy gradients coincide with those in~\cite{zhao2023global}. 
Overall, the theoretical characterization parallels that of the standard LQR problem~\cite{yang2019provably} and the mean-field linear–quadratic control setting~\cite{fu2019actor}.
\end{remark}

\begin{remark}
To the best of our understanding, the derivation of the action-value function in~\eqref{Q:form} and its analytical connection to the policy gradient for the RC-LQR problem has not been explicitly presented in prior work. 
Lemma~\ref{lem:policy_gradient} forms an important step toward developing the proposed actor–critic framework.
\end{remark}

\subsection{Multi-Time-Scale Actor-Critic Algorithm}

Next, we proceed to design an actor–critic algorithm for problem~\eqref{duality} within a primal–dual framework. 
Specifically, the primal variable $X$ is updated along the negative gradient direction of the Lagrangian $\mathcal{L}(X,\mu)$, while the dual variable $\mu$ is updated along its gradient ascent direction. 
In the model-free setting, the exact policy gradients are unavailable. 
To address this challenge, we estimate the action-value function $Q^{\mu}_{X}$ in the \emph{critic step}. 
Then, using Lemma~\ref{lem:policy_gradient}, we compute the policy gradient of $\mathcal{L}(X,\mu)$ with respect to $X$ and update the primal variable in the \emph{actor step}. 
Finally, the dual variable $\mu$ is updated based on the estimated constraint function $\mathcal{J}_{c}(X)$.

For a given policy $X \in \mathcal S$ and a multiplier $\mu \in \mathbb R_+$, we learn $Q_X^\mu$ via a temporal-difference (TD) algorithm as follows.

For a pair $(x, u)$, we define the feature vector as
\begin{equation*}
\begin{aligned}
\label{feature}
\psi(x, u) \!=\! \begin{bmatrix}\!
\phi(x, u) \\
2x \\
2u \!
\end{bmatrix}\!\!,
~{\rm where}~
\phi(x, u) \!=\! {\rm svec} 
\begin{bmatrix}
\begin{pmatrix}
x \\
u
\end{pmatrix}\!
\begin{pmatrix}
x  \\
u 
\end{pmatrix}^\top
\end{bmatrix}\!.
\end{aligned}
\end{equation*}
In addition, we denote
\begin{equation}
\begin{aligned}
\label{para:Q}
\Theta_X^{\mu} =
\begin{bmatrix}
{\rm svec} 
(\Upsilon_{K,\mu}) \\
p_{X,\mu} \\
q_{X,\mu}
\end{bmatrix}.
\end{aligned}
\end{equation}
Recalling (\ref{Q:form}) gives
\begin{equation}
\label{Q:reform}
Q^{\mu}_X = \psi(x, u)^\top \Theta_X^{\mu} + \varsigma_{X,\mu}.
\end{equation}
where $ \varsigma_{X,\mu}$ is a scalar independent of $x$ and $u$.
Clearly, to estimate $\Upsilon_{K,\mu}$ and $q_{X,\mu}$, it suffices to compute $\Theta_X^{\mu}$.

Note that the action-value function $Q^{\mu}_X$ satisfies a Bellman equation as
\begin{equation*}
\begin{aligned}
Q^{\mu}_X(x_t, u_t)=&~ c_\mu(x_t, u_t)-\mathcal L(X, \mu)\\
&+ \mathbb E \big[Q^{\mu}_X(x_{t+1}, u_{t+1})\mid x_t, u_t \big].
\end{aligned}
\end{equation*}
Substituting (\ref{Q:reform}) into the above equation, we derive
\begin{equation}
\begin{aligned}
\label{Bell:Q}
\psi(x_t, u_t)^\top \Theta^{\mu}_X =&~ c_\mu(x_t, u_t) - \mathcal L(X, \mu) \\
&+  \mathbb E \big[\psi(x_{t+1}, u_{t+1})^\top \cdot \Theta^{\mu}_X \mid x_t, u_t \big].
\end{aligned}
\end{equation}

For notational convenience, we denote
$c_{\mu,t} := c_{\mu}(x_t, u_t)$, 
$\psi_t := \psi(x_t, u_t)$, 
and $\mathcal{L}_t := \mathcal{L}(X_t, \mu_t)$. 
Let $\hat{\mathcal{L}}_t$ be an estimate that tracks the long-term Lagrangian value $\mathcal{L}(X,\mu)$, and let $\vartheta^{\mu}_{X,t}$ denote the estimate of the parameter vector $\Theta^{\mu}_{X}$ at iteration~$t$, where
\begin{equation}
\label{Est:Psi}
\vartheta^{\mu}_{X,t}
=
\big[\big(\mathrm{svec}(\widehat{\Upsilon}_{(K,\mu),t})\big)^{\!\top},~
\hat{p}^{\top}_{(X,\mu),t},~
\hat{q}^{\top}_{(X,\mu),t}\big]^{\!\top}.
\end{equation}
Here, $\widehat{\Upsilon}_{(K,\mu),t}$, $\hat{p}_{(X,\mu),t}$, and $\hat{q}_{(X,\mu),t}$ denote the estimates of $\Upsilon_{K,\mu}$, $p_{X,\mu}$, and $q_{X,\mu}$, respectively.

Then the TD algorithm is designed as
\begin{align}
\hat{\mathcal L}_{t+1} &= (1 - \alpha_t)\hat{\mathcal L}_t + \alpha_t \cdot c_{\mu, t}, \label{TD0:cost} \\
\vartheta^{\mu}_{X, t+1} &= \vartheta^{\mu}_{X, t} + \alpha_t \cdot \delta_t  \cdot \psi_t, \label{TD0:Q}
\end{align}
where $\alpha_t > 0$ is the stepsize, and the TD error $\delta_t$ is
\begin{equation}
\label{delta}
\delta_t =  c_{\mu, t} - \hat{\mathcal L}_t 
+ (\psi_{t+1} - \psi_t)^\top \vartheta^{\mu}_{X, t}.
\end{equation}

The TD algorithm follows the general structure of those in~\cite{sutton2018reinforcement, bhatnagar2009natural, zhang2018fully}, 
but is tailored to the RC-LQR problem, where both the feature vector and the action-value function are parameterized in problem-specific forms.

According to~\eqref{policy:gradient}, computing the policy gradient requires an estimate of $\Phi_{X}$. 
Analogous to the evolution of $\hat{\mathcal{L}}_t$, we estimate $\Phi_{X}$ via
\begin{equation}
\label{est:Phi}
\hat{\Phi}_{X,t+1}
= (1 - \alpha_t)\hat{\Phi}_{X,t}
+ \alpha_t\, z_t z_t^{\top},
\end{equation}
where $\hat{\Phi}_{X}$ tracks $\Phi_{X}$ in~\eqref{Phi:eq}, and $z_t = [x_t^{\top}, -1]^{\top}$.

It is worth noting that our focus is on computing the exact policy gradient for the update of $X$, and hence~\eqref{est:Phi} is used to estimate $\Phi_{X}$. 
In contrast, if a natural policy gradient algorithm were adopted, as in~\cite{yang2019provably, kakade2001natural, bhatnagar2009natural}, such estimation would not be required.

Recalling (\ref{lagrangian}) gives $\nabla_{\mu} \mathcal L(X, \mu) = \mathcal J_c(X) - \bar \iota$.
For the evolution of $\mu$, we estimate $\mathcal J_c(X)$ by
\begin{equation}
\label{est:const}
\hat{\mathcal J}^c_{X, t+1} = (1 - \alpha_t)\hat{\mathcal J}^c_{X, t} + \alpha_t \cdot o_t,
\end{equation}
where $o_t = x_t^\top QWQ x_t + 4x_t^\top QM_3- \bar\iota$.

For the actor step, we first recover $\hat \Upsilon_{(K,\mu), t}$ and $\hat q_{(X,\mu), t}$ from $\vartheta_{X, t}^{\mu}$ according to (\ref{Est:Psi}).
Then the estimates of $E_{X, \mu}$ and $G_{X, \mu}$, denoted by $\hat E_{(K, \mu), t}$ and $\hat G_{(K, \mu), t}$, are computed via
\begin{equation*}
\begin{aligned}
\hat E_{(K, \mu), t} &= \hat\Upsilon_{(K,\mu),t}^{22} K_t - \hat \Upsilon_{(K,\mu),t}^{21} \\
\hat G_{(X, \mu), t} &= \hat\Upsilon_{(K,\mu),t}^{22} b_t - \hat q_{(X, \mu), t}.
\end{aligned}
\end{equation*}
Define $\hat H_{(X, \mu), t} = [\hat E_{(K, \mu), t}, \hat G_{(X, \mu), t}]$.
Based on these estimates, the policy gradient algorithm for $X$ is designed as
\begin{equation}
\begin{aligned}
\label{policy-gradient}
X_{t+1} = \Gamma_{\mathcal X}\big[X_t - \beta_t \hat H_{(X, \mu), t} \hat \Phi_{X, t} \big],
\end{aligned} 
\end{equation}
where $\beta_t > 0$ is the stepsize, and $\Gamma_{\mathcal X}(\cdot)$ denotes the Euclidean projection operator to ensure the boundedness of the sequence $\{X_t\}$. The set $\mathcal X$ is compact and convex, satisfying 
$$\mathcal X \subset \{X \mid \rho(A-BK) \le 1 - \epsilon\},$$
for some $\epsilon \in (0, 1)$, and contains  the optimal policy $X^*$.

In fact, projection operators are widely used to stabilize the stochastic approximation algorithms \cite{kushner2003stochastic}. The aforementioned requirements are standard in the analysis of actor-critic algorithms \cite{bhatnagar2009natural, zhang2018fully, fu2019actor}.
Notably, the projection is mainly introduced for theoretical analysis and may not be necessary in practical implementations.

Finally, the dual variable $\mu$ is updated using a projected gradient ascent scheme:
\begin{equation}
\label{dual:ascent}
\mu_{t+1}
= \Gamma_{\mathbb{R}_{+}}\!\left[\mu_t + \gamma_t \big(\hat{\mathcal{J}}^{c}_{X,t} - \bar{\iota}\big)\right],
\end{equation}
where $\gamma_t > 0$ is the stepsize, and $\Gamma_{\mathbb{R}_{+}}(\cdot)$ denotes the Euclidean projection operator onto the nonnegative orthant, ensuring $\mu_t \ge 0$ for all~$t$.

We impose the following standard conditions on the stepsizes.

\begin{assumption}
\label{ass:stepsize}
Let $\{\alpha_t\}$, $\{\beta_t\}$, and $\{\gamma_t\}$ be positive stepsize sequences satisfying
\begin{equation}
\begin{aligned}
\label{stepsize}
&\sum\nolimits_t \alpha_t
= \sum\nolimits_t \beta_t
= \sum\nolimits_t \gamma_t
= \infty, \\
&\sum\nolimits_t (\alpha_t^2 + \beta_t^2 + \gamma_t^2)
< \infty, \\
&\lim\nolimits_{t \to \infty} \frac{\beta_t}{\alpha_t} = 0, 
\quad
\lim\nolimits_{t \to \infty} \frac{\gamma_t}{\beta_t} = 0.
\end{aligned}
\end{equation}
\end{assumption}

\begin{remark}
The first two conditions in Assumption~\ref{ass:stepsize} are standard in stochastic approximation theory~\cite{kushner2003stochastic}. 
The final condition enforces a separation of time scales: the critic updates occur at the fastest rate, followed by the actor, while the dual ascent step evolves the slowest. 
A similar assumption has been adopted in~\cite{borkar2005actor, bhatnagar2010actor}.
\end{remark}

We summarize the proposed multi–time–scale actor–critic algorithm in Algorithm~\ref{alg}. 
This framework extends the conventional two–time–scale actor–critic schemes~\cite{bhatnagar2009natural, zhang2018fully, zeng2024two} by incorporating an additional time scale associated with the dual variable. 
Unlike the bilevel optimization structures considered in~\cite{yang2019provably, fu2019actor}, all decision variables are updated concurrently. 
Furthermore, the method differs from the sample–based policy gradient algorithm in~\cite{zhao2023global} in that neither the global cost nor the constraint function needs to be explicitly sampled.

\begin{algorithm}[t]
\caption{Actor-Critic Algorithm for RC-LQR} 
\label{alg}
\begin{algorithmic}
\STATE \textbf{Input:}

\STATE \hspace{0.2cm} The initial policy $X_0$, multiplier $\mu_0$, and state $x_0$.

\STATE \hspace{0.2cm} The initial values of $\hat{\mathcal L}_0$, $\hat \Phi_{X, 0}$, $\hat{\mathcal J}^c_{X, 0}$ and $\vartheta^{\mu}_{X, 0}$.

\STATE \hspace{0.2cm} The stepsize sequences $\{\alpha_t\}$, $\{\beta_t\}$ and $\{\gamma_t\}$.

\STATE {\textbf{Repeat:}}

\STATE \hspace{0.2cm} Execute $u_t = -K_t x_t + b_t + \sigma \eta_t$.

\STATE \hspace{0.2cm} Observer the state $x_{t+1}$ and cost $c_{\mu, t}$.

\STATE \hspace{0.2cm} Update $\hat{\mathcal L}_{t+1} = (1 - \alpha_t)\hat{\mathcal L}_t + \alpha_t \cdot c_{\mu, t}$.

\STATE \hspace{0.2cm} Update $\hat \Phi_{X, t+1} = (1 - \alpha_t)\hat \Phi_{X, t} +  \alpha_t \cdot z_t z_t^\top$.

\STATE \hspace{0.2cm} Update $\hat{\mathcal J}^c_{X, t+1} = (1 - \alpha_t)\hat{\mathcal J}^c_{X, t} + \alpha_t \cdot o_t$.

\STATE \hspace{0.2cm} Compute the TD error
$\delta_t =  c_{\mu, t} - \hat{\mathcal L}_t 
+ (\psi_{t+1} - \psi_t)^\top \vartheta^{\mu}_{X, t}$.

\STATE \hspace{0.2cm} 
\textbf{Critic step:}
$\vartheta^{\mu}_{X, t+1} = \vartheta^{\mu}_{X, t} + \alpha_t \cdot \delta_t  \cdot \psi_t$.

\STATE \hspace{0.2cm} 
\textbf{Actor step:}
\begin{equation*}
\begin{aligned}
X_{t+1} = \Gamma_{\mathcal X}\big[X_t - \beta_t  \hat H_{(X, \mu), t}  \hat \Phi_{X, t}\big].
\end{aligned} 
\end{equation*}

\STATE \hspace{0.2cm} \textbf{Dual ascent:}
$\mu_{t+1} = \Gamma_{\mathbb R_+} \big[\mu_t + \gamma_t (\hat{\mathcal J}^c_{X, t} - \bar \iota) \big]$.

\STATE \textbf{Output:}~Policy $X$ and the multiplier $\mu$.
\end{algorithmic}
\end{algorithm}

\section{THEORETICAL RESULTS}

In this section, we analyze the convergence properties of Algorithm~\ref{alg} using the ordinary differential equation (ODE) method developed in~\cite{kushner2003stochastic, borkar2008stochastic}. 
The analysis is based on the multi–time–scale stochastic approximation framework of~\cite{borkar1997stochastic}. 
In particular, under Assumption~\ref{ass:stepsize}, the results in~\cite{borkar1997stochastic} ensure the convergence of the critic update when the policy $X$ and multiplier $\mu$ are held fixed. 
Subsequently, we establish the convergence of the policy update for a fixed $\mu$, and finally analyze the convergence behavior of the dual variable~$\mu$.

For a fixed policy $X$ and multiplier $\mu$, we first analyze the convergence of the critic step as follows.
\begin{theorem}
\label{thm:TD}
Let Assumptions \ref{ass:mat}-\ref{ass:stepsize} hold.
Consider the sequences $\{\hat {\mathcal L}_t\}$, $\{\vartheta^{\mu}_{X,t}\}$, $\{\hat \Phi_{X, t}\}$ and $\{\hat{\mathcal J}^c_{X, t}\}$ generated by (\ref{TD0:cost}), (\ref{TD0:Q}), (\ref{est:Phi}) and (\ref{est:const}), respectively.
For any fixed policy $X \in \mathcal{X}$ and multiplier $\mu \in \mathbb{R}_+$, it follows that
\begin{equation*}
\begin{aligned}
&\hat {\mathcal L}_t \rightarrow \mathcal L(X, \mu), ~
\vartheta^{\mu}_{X,t} \rightarrow \vartheta^{\mu, *}_X, ~\\
&\hat \Phi_{X, t} \rightarrow \Phi_X, 
~{\rm and}~\hat{\mathcal J}^c_{X, t} \rightarrow \mathcal J_c(X) ~{\rm a.s.}
\end{aligned}
\end{equation*}
where $\mathcal L(X, \mu)$, $\Phi_X$ and $\mathcal J_c(X)$ are defined in (\ref{lagrangian}), (\ref{Phi:eq}) and (\ref{reform}), $\vartheta^{\mu, *}_X$ is the unique solution to
\begin{equation*}
\begin{aligned}
\mathbb E_{(x_t, u_t)}\big[\psi_t(\psi_t - \psi_{t+1})^\top\big] 
\cdot \vartheta^{\mu, *}_X& \\
+ \mathbb E_{(x_t, u_t)} [\psi_t] \cdot \mathcal L(X, \mu)
&= \mathbb E_{(x_t, u_t)}[c_{\mu,t} \psi_t].
\end{aligned}
\end{equation*}
and moreover, $\mathbb E_{(x_t, u_t)}[\cdot]$ denotes the expectation with respect to $x_t \sim \nu_X$ and $u_t \sim \pi_X^{\mu}(\cdot \mid x_t)$.
\end{theorem}

\emph{Proof:}
The local update in (\ref{TD0:cost}) forms a stochastic approximation iteration, whose asymptotic behavior is captured by the ODE
\begin{equation*}
\label{ode:cost}
\dot {\hat{\mathcal L}} = - \hat{\mathcal L} + \mathcal L(X, \mu).
\end{equation*}

Substituting $\delta_t$ from (\ref{delta}) into (\ref{TD0:Q}), we derive
$\vartheta^{\mu}_{X, t+1} = \vartheta^{\mu}_{X, t} + \alpha_t \cdot \big[(c_{\mu, t} - \hat{\mathcal L}) \psi_t 
+ \psi_t (\psi_{t+1} - \psi_t)^\top \vartheta^{\mu}_{X, t} \big].
$
Then the ODE characterizing the asymptotic behavior of this recursion is given by
\begin{equation*}
\label{ODE:Q}
\dot \vartheta^{\mu}_X = \Xi \cdot \vartheta^{\mu}_X - \bar \psi \hat{\mathcal L} + l,
\end{equation*}
where $\Xi = \mathbb E_{(x_t, u_t)} \big[\psi_t (\psi_{t + 1} - \psi_t)^\top\big]$, $\bar \psi = \mathbb E_{(x_t, u_t)}[\psi_t]$, and 
$l = \mathbb E_{(x_t, u_t)}[c_{\mu, t}\psi_t]$.
Denote by $\omega = [\vartheta^{\mu, \top}_X, \hat{\mathcal L}]^\top$, $\hat \Xi = [\Xi, -\bar \psi; ~0, -1]$ and $\hat l = [l^\top, \mathcal L(X, \mu)]^\top$.
It follows that
\begin{equation}
\label{TD:ODE}
\dot \omega = \hat \Xi \cdot \omega + \hat l.
\end{equation}

Let $f(\omega)$ be the right-hand-side of (\ref{TD:ODE}). Clearly, $f(\omega)$ is Lipschitz continuous in $\omega$.
Define $f_{\infty}(\omega) = \lim_{r\to \infty} \frac {f(r\omega)}{r}$.
It follows that
$f_{\infty}(\omega) = \hat \Xi \cdot \omega$.
We next analyze the performance of a dynamics given by
\begin{equation}
\label{TD:ODE:inf}
\dot \omega = \hat \Xi \cdot \omega.
\end{equation}

It follows directly that
\[
\Xi + \Xi^\top 
= -\,\mathbb{E}_{(x_t, u_t)} \!\left[(\psi_{t+1} - \psi_t)(\psi_{t+1} - \psi_t)^\top\right]
\preceq 0,
\]
which implies that $\Xi$ is negative semidefinite.

Next, we establish that $(\Xi + \Xi^\top) \prec 0$. 
Let $\mathcal{F}^1_t$ denote the filtration 
$\mathcal{F}^1_t = \sigma(x_0, w_0, \ldots, w_{t-1}, \eta_0, \ldots, \eta_t)$,
which forms an increasing $\sigma$-algebra with respect to time~$t$. 
Define $\Delta \psi_t = \psi_{t+1} - \psi_t$.
For any nonzero vector $y \in \mathbb{R}^{\dim(\Delta \psi_t)}$, 
let $r_y = y^\top \Delta \psi_t$. 
Then, it holds that
\begin{equation*}
\begin{aligned}
-\,y^\top (\Xi + \Xi^\top) y 
&= \mathbb{E}_{(x_t, u_t)} [r_y^2] 
= \mathbb{E}_{(x_t, u_t)} \!\left[\mathbb{E}[r_y^2 \mid \mathcal{F}^1_t]\right].
\end{aligned}
\end{equation*}

Since $\mathbb E \left[(w_t - \bar w)(w_t - \bar w)^\top\right] \succ 0$, it follows from (\ref{LTI:sys}) that for any $y_x \in \mathbb R^n$, $\mathbb E \big[\Vert y_x^\top(x_{t+1} - x_t) \Vert^2 \mid \mathcal F^1_t \big] = 0$ only if $y_x = 0$.
Using a similar argument, we have
$\mathbb E[r_y^2 \mid \mathcal F^1_t] = 0$ if and only if $y = 0$ since $u_t = -Kx_t + b + \sigma \eta_t$, where $\eta_t \sim \mathcal N(0, I_m)$ and $\sigma > 0$.
Consequently, $-\,y^\top (\Xi + \Xi^\top) y > 0$ for all $y \neq 0$, 
implying that $(\Xi + \Xi^\top) \prec 0$. 
Therefore, $\Xi$ is a Hurwitz matrix.

By definition, the eigenvalues of $\hat{\Xi}$ consist of those of $\Xi$ together with an additional eigenvalue at $-1$. 
Since all eigenvalues of $\Xi$ have negative real parts, $\hat{\Xi}$ is Hurwitz. 
Hence, the ODE in~\eqref{TD:ODE:inf} admits the origin as its unique globally asymptotically stable equilibrium.

Let $\mathcal F^2_t$ be the $\sigma$-algebra given by
$\mathcal F^2_t = \sigma(\hat{\mathcal L}_s, \vartheta^{\mu}_{X, s}, s \le t)$.
Define 
$M^1_{t+1} = c_{\mu, t} - \mathbb E[c_{\mu, t} \mid \mathcal F^2_t]$, and
$M^2_{t+1} = \delta_t  \cdot \psi_t - \mathbb E [\delta_t  \cdot \psi_t \mid \mathcal F^2_t]$. It follows that
both $\{M^1_t, \mathcal F^2_t: t \ge 0\}$ and $\{M^2_t, \mathcal F^2_t: t \ge 0\}$ are martingale difference sequences.
For any $X \in \mathcal X$ and $\mu \ge 0$, it is easy to see that
\begin{equation*}
\begin{aligned}
\mathbb E \big[\Vert M_{t+1}^1 \Vert^2 \mid \mathcal F^2_t \big] \le C_1 \cdot \big(1 + {\hat{\mathcal L}}_t^2 + \Vert \vartheta^\mu_{X, t} \Vert^2\big), \\
\mathbb E \big[\Vert M_{t+1}^2\Vert^2 \mid \mathcal F^2_t \big] \le C_2 \cdot \big(1 + {\hat{\mathcal L}}_t^2 + \Vert \vartheta^\mu_{X, t} \Vert^2\big),
\end{aligned}
\end{equation*}
for some $C^1, C^2 < \infty$.

Recalling Lemma~\ref{lem:SA:bound} yields $\sup_t \Vert \omega_t \Vert < \infty$ a.s. 
Then, by Lemma~\ref{lem:SA:con}, we have $\omega_t \!\to\! \omega^*$ a.s., 
where $\omega^* = (\vartheta^{\mu,*}_{X}, \hat{\mathcal{L}}^*)$ denotes the equilibrium of~\eqref{TD:ODE} satisfying
\begin{equation*}
\Xi\,\vartheta^{\mu,*}_X - \bar{\psi}\,\hat{\mathcal{L}}^* + l = 0,
\qquad
\hat{\mathcal{L}}^* = \mathcal{L}(X, \mu).
\end{equation*}
Moreover, the equilibrium $\vartheta^{\mu,*}_X$ is unique since $\Xi$ is Hurwitz.

By an analogous argument, we can show that 
$\hat{\Phi}_{X,t} \!\to\! \Phi_X$ and 
$\hat{\mathcal{J}}^{c}_{X,t} \!\to\! \mathcal{J}_c(X)$ a.s. 
This completes the proof.
$\hfill\square$

\begin{remark}
\label{rmk:TD}
The quantity $\vartheta^{\mu,*}_X$ represents the exact solution of the Bellman equation~\eqref{Bell:Q} corresponding to the given policy $X$ and multiplier~$\mu$. 
Theorem~\ref{thm:TD} therefore shows that the proposed TD algorithm accurately estimates the associated action–value function.
\end{remark}

Next, we turn to the convergence analysis of the policy sequence $\{X_t\}$ in Algorithm~\ref{alg}, assuming a fixed multiplier~$\mu$.

\begin{theorem}
\label{thm:policy}
Suppose Assumptions~\ref{ass:mat}–\ref{ass:stepsize} hold.
Consider the sequence $\{X_t\}$ generated by~\eqref{policy-gradient}.
For any fixed $\mu \in \mathbb{R}_+$, the iterates $\{X_t\}$ converge a.s. to an equilibrium point of
\begin{equation}
\label{actor:dyn}
\dot{X} = \widehat{\Gamma}_{\mathcal{X}}\big[-h_{\mu}(X)\big],
\end{equation}
where $h_{\mu}(X) = H_{X,\mu}\Phi_X$, and moreover, both $H_{X,\mu}$ and $\Phi_X$ are defined in~\eqref{policy:gradient}.
\end{theorem}

\emph{Proof:}
Let $\mathcal F^3_t$ be the $\sigma$-field generated by $\{X_s, s \le t\}$.
It is clear that the recursion (\ref{policy-gradient}) reads as
\begin{equation}
\begin{aligned}
\label{actor:ite}
X_{t+1} &= \Gamma_{\mathcal X} \big[X_t - \beta_t h_\mu(X)  + \beta_t \xi^1_t\big],
\end{aligned}
\end{equation}
where $\xi^1_t = h_\mu(X) - \hat H_{(X, \mu), t} \hat \Phi_{X, t}$.
Recalling Theorem \ref{thm:TD} and Remark \ref{rmk:TD}, gives $\lim_{t \to \infty}\xi_t^1 = 0$. By Lemma \ref{lem:SA:proj}, the iteration (\ref{actor:ite}) can be captured by the ODE 
$\dot X = \widehat \Gamma_{\mathcal X}[-h_{\mu}(X)].$

Construct a Lyapunov function candidate as
$
S_1(X) = \mathcal L(X, \mu) - \mathcal L(X^*, \mu).
$
For a fixed $\mu$, we have
\begin{equation*}
\begin{aligned}
\dot S_1 &= \big\langle h_{\mu}(X), \widehat \Gamma_{\mathcal X}[-h_{\mu}(X) \big\rangle
= - \big\Vert \widehat \Gamma_{\mathcal X}[-h_{\mu}(X)] \big\Vert^2 \le 0.
\end{aligned}
\end{equation*}
where the second equality follows from the Moreau decomposition.
As a consequence, $X_t$ approaches the equilibrium of (\ref{actor:dyn}) a.s. by Lemma \ref{lem:SA:proj}, and this completes the proof.
$\hfill\square$

\begin{remark}
Since the optimal policy $X^*$ lies in $\mathcal{X}$, the equilibrium set of~\eqref{actor:dyn} contains $X^*$. 
We note that Theorem~\ref{thm:policy} is consistent with the convergence results established in~\cite{bhatnagar2009natural, zhang2018fully}, and provides the standard convergence guarantee expected for actor–critic algorithms.
\end{remark}

Next, we analyze the convergence of the multiplier sequence $\{\mu_t\}$ in Algorithm~\ref{alg}.

\begin{theorem}
\label{thm:dual}
Suppose Assumptions~\ref{ass:mat}–\ref{ass:stepsize} hold.  
Then, the sequence $\{\mu_t\}$ generated by~\eqref{dual:ascent} converges a.s. to $\mu^*$,  
where $\mu^* \in \arg\max_{\mu \ge 0} \mathcal{L}(\tilde{X}_{\mu^*}, \mu)$,  
and $\tilde{X}_{\mu^*}$ denotes an equilibrium of~\eqref{actor:dyn} corresponding to~$\mu^*$.
\end{theorem}

\emph{Proof:}
By Assumption~\ref{ass:stepsize}, the recursion~\eqref{policy-gradient} evolves on a faster timescale than~\eqref{dual:ascent}. 
Moreover, for any fixed $\mu \in \mathbb{R}_+$, Theorem~\ref{thm:policy} ensures that $X_t$ converges a.s. to an equilibrium point $\tilde{X}_{\mu}$ of~\eqref{actor:dyn}. 
Substituting this into~\eqref{dual:ascent} yields
\begin{equation}
\label{ite:dual}
\mu_{t+1}
= \Gamma_{\mathbb{R}_+}\!\left[\mu_t + \gamma_t \big(\mathcal{J}_c(\tilde{X}_{\mu}) - \bar{\iota} + \xi^2_t\big)\right],
\end{equation}
where $\xi^2_t = \hat{\mathcal{J}}^c_{X_{\mu},t} - \mathcal{J}_c(\tilde{X}_{\mu})$, and $\lim_{t \to \infty} \xi^2_t = 0$.

In light of Lemma \ref{lem:SA:proj}, the ODE associated with (\ref{ite:dual}) is 
$$
\dot \mu = \widehat \Gamma_{\mathbb R_+} \big[\mathcal J_c(\tilde X_{\mu}) - \bar \iota \big].
$$

Construct a Lyapunov function candidate as
$$
S_2 = \frac 12 \Vert \mu - \mu^*\Vert^2.
$$
where $\mu^* \in \arg\max \mathcal L(\tilde X_{\mu}, \mu)$. Then 
\begin{equation*}
\begin{aligned}
\dot S_2 \le \langle \mu - \mu^*, \mathcal J_c(\tilde X_{\mu}) - \bar \iota\rangle
= \mathcal L(\tilde X_{\mu}, \mu) - \mathcal L(\tilde X_{\mu}, \mu^*) \le 0,
\end{aligned}
\end{equation*}
where the first inequality follows from the properties of the projection operator. 
Moreover, $\dot{S}_2 = 0$ only if 
$\mu \in \arg\max_{\mu \ge 0} \mathcal{L}(\tilde{X}_{\mu}, \mu)$. 
Hence, $\mu_t \!\to\! \arg\max_{\mu \ge 0} \mathcal{L}(\tilde{X}_{\mu}, \mu)$ a.s. 
This completes the proof.
$\hfill\square$

\begin{remark}
In contrast to the two–time–scale actor–critic algorithms in~\cite{bhatnagar2009natural, zhang2018fully, zeng2024two}, 
the proposed method incorporates an additional dual variable to enforce the constraint and adopts a multi–time–scale learning structure. 
Moreover, we establish asymptotic convergence of the overall algorithm, thereby extending the convergence guarantees of the two–time–scale framework.
\end{remark}

\section{SIMULATION}

This section provides numerical simulations to verify the performance of Algorithm \ref{alg}.

Consider the dynamical system (\ref{LTI:sys}) with
\begin{equation*}
\begin{aligned}
A = \begin{bmatrix}
1 &1 &0 &0 \\
0 &1 &0 &0 \\
0 &0 &1 &1 \\
0 &0 &0 &1
\end{bmatrix}\!, ~{\rm and}~
B = \begin{bmatrix}
1 &0 \\
1 &0 \\
0 &1 \\
0 &1
\end{bmatrix}\!.
\end{aligned}
\end{equation*}
The quadratic cost function (\ref{cost}) is defined with
$Q = {\rm diag}(0.5, 0.1, 0.1, 0.5)$ and $R = 0.2 I_2$.
The noise sequence $\{w_t\}$ is specified as $w_t = B(\tilde w_t + \hat w_t)$, 
where $\tilde w_t = [\tilde w_t^1, \tilde w_t^2]^\top$, $\hat w_t = [\hat w_t^1, \hat w_t^2]^\top$,
$\tilde w_t^1$ follows a mixed Gaussian distribution composed of $\mathcal N(5, 8)$ and $\mathcal N(8, 10)$ with weights $0.3$ and $0.7$, $\tilde w_t^2$ follows a Gaussian distribution as $\mathcal N(0, 4)$, both $\hat w^1_t$ and $\hat w^2_t$ are from a uniform distribution as ${\rm Unif}(0, 0.5)$, and $\iota = 40$.

\begin{figure}[htp]
\centering
\includegraphics[scale=0.31]{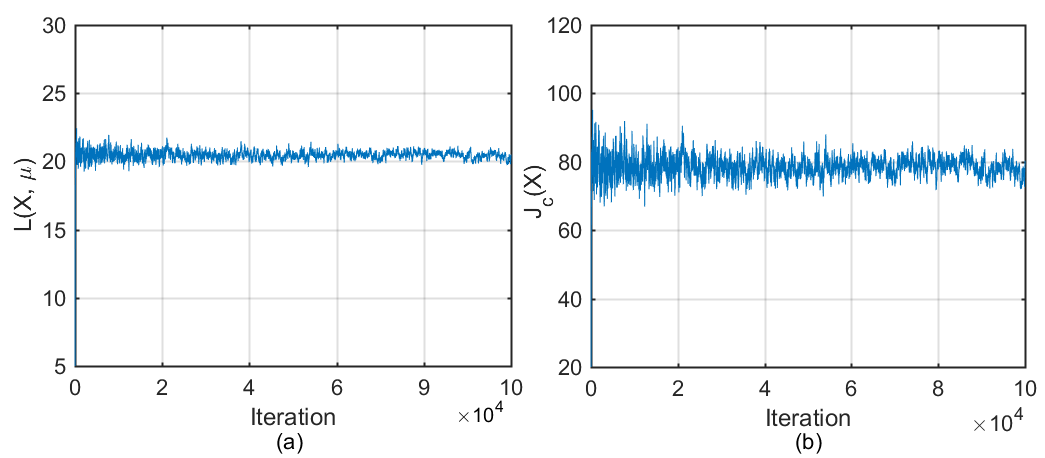}
\vspace{-0.2cm}
\caption{(a) The trajectory of $\hat{\mathcal L}_t$. (b) The trajectory of $\hat{\mathcal J}^c_{X, t}$.}
\label{Fig1}
\end{figure}

Fig. \ref{Fig1} shows the trajectories of $\hat{\mathcal L}_t$ and $\hat{\mathcal J}^c_{X, t}$. It is clear that 
both of them converges as proved in Theorem \ref{thm:TD}.
Fig. 2 shows the trajectories of $\Vert X_t - X^*\Vert_F$ and the multiplier $\mathcal \mu_t$, where $X^*$ is computed by Algorithm in \cite{zhao2021infinite}.
Thus, $X_t$ approaches the optimal solution, and $\mu_t$ is also convergent.

\begin{figure}[htp]
\centering
\includegraphics[scale=0.31]{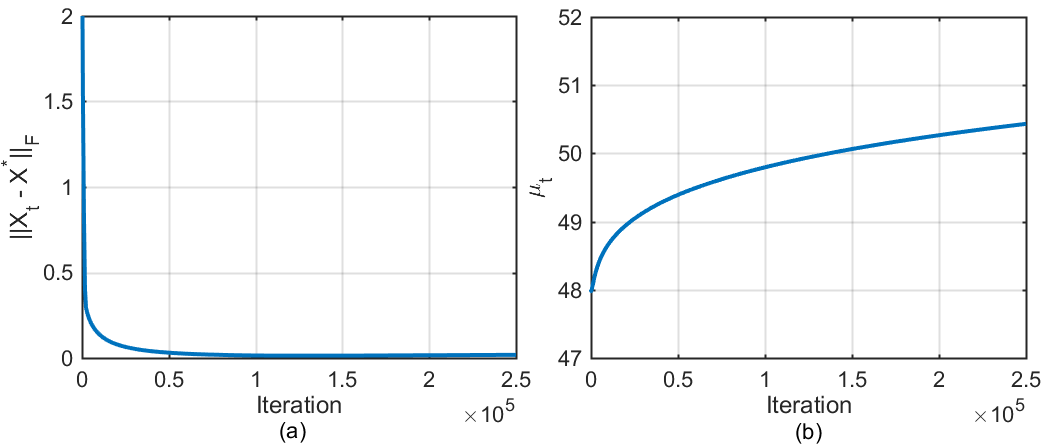}
\vspace{-0.2cm}
\caption{(a) The trajectory of $\Vert X_t - X^*\Vert_F$. (b) The trajectory of $\mu_t$.}
\label{Fig2}
\end{figure}

\section{CONCLUSION}

In this paper, we investigated an RC–LQR problem with unknown system dynamics. 
Unlike the standard LQR formulation, we imposed a constraint on the one–step predictive variance of the state to enhance safety in practical applications. 
To address this problem, we developed a primal–dual actor–critic algorithm, in which the critic employed a temporal–difference method to estimate the action–value function, the actor updated the policy via a policy gradient step, and the dual variable was adjusted through a gradient ascent update. 
We designed multi–time–scale learning rates to ensure convergence of the learning process, and established asymptotic convergence using stochastic approximation theory. 
Finally, we demonstrated the effectiveness of the proposed approach through numerical simulations.

\section*{References}
\vspace{-1.5em}
\bibliographystyle{IEEEtran}
\bibliography{references,IDS_Publications_10132025}

\end{document}